\documentclass{article}
\usepackage{amsmath}
\usepackage{bbm}
\usepackage{hyperref}
\usepackage{sectsty} 
\usepackage{enumitem}
\usepackage{multicol}
\sectionfont{\normalfont\normalsize\centering\scshape\underline}
\title{A Cubic Composite Test} 
\author{Pierre Laurent, Paul Underwood}

\begin{document}
\maketitle
\date

\begin{abstract}
  A single parameter cubic composite test for odd positive integers is given which relies on the discriminant always being a square integer.
  This test has no known counterexample despite extensive verifications. As well as a comparison with
  the Baillie-PSW tests,
  a related quadratic composite test is briefly examined which also has no known counterexample.
\end{abstract}

\section{Introduction}

Grantham gives an excellent introduction in his paper Frobenius Pseudoprimes \cite{GRANTHAM},
in which various pseudoprime tests are performed over polynomials in one variable.
Therein simultaneous modulo $n$ and variously modulo a polynomial are considered.
We formally define working \!\!\!\! $\pmod{n,f_a}$ as working out our arithmetic in the quotient ring
$\mathbbm{Z}_n[x]/f_a$ where $f_a$ is a polynomial in $x$ with specific integer coefficients.
Note that both $n\equiv 0 \pmod{n,f_a}$ and $f_a\equiv 0 \pmod{n,f_a}$ $(*)$.

In some sense we herein build upon the pseudoprimes paper by Adams and Shanks \cite{SHANKS}.
We shall make a computational cost comparison of the cubic test with the Baillie-PSW tests.
We shall also take a brief look at a quadratic composite test based on the
charateristic equation $x^2-2x-4=0$.

This paper has sections on each of the conditions to construct a cubic composite test algorithm.
It does not consider sieving for primes, trial division and other techniques for finding factors,
nor other quick tests, all of which might speed up batch testing of candidate composites.

\section{First g.c.d.}
We consider depressed cubic polynomials of the form $f_a = x^3-ax-a$ simply with $\gcd(a,n)=1$.

\section{Second g.c.d.}
We henceforth restrict $f_a$ with the parametric equation $a=7+k(k-1)$.
The cubic polynomials $x^3-ax-a$ have discriminants $4a^3-27a^2$.
Since squares can be factored out of Kronecker symbols, only $4a-27$ can be considered.
Making the substitution in terms of $k$, this expression becomes $28+4k(k-1)-27$ which
is equal to $4k(k-1)+1$; is equal to $(2k-1)^2$.
Thus the Kronecker symbols of the discriminants of $f_a$ are never negative.
This implies that each $f_a$ has three real roots.
In practice the greatest common divisors can be taken instead: $\gcd(2k-1,n)$.
Only a g.c.d. of $1$ will be of interest to us for primality testing purposes.

\section{Third g.c.d.}

Consider the following necessary characteristic polynomial for the companion matrix of $x^3-ax-a$:
\[
\left|
z-
\frac
{\left(
\begin{array}{ccc}
  0& a& a\\
  1& 0& 0\\
  0& 1& 0
\end{array}
\right)^2}{-a}
\right|
=z^3+2z^2+z+\frac{1}{a}.
\]
If $a=\frac{1}{2}$ then the right hand side of the equation factors as $(z+2)(z^2+1)$.
So $\gcd(2a-1,n)=1$ can be checked.
This g.c.d. test arose historically while developing the cubic test of this
paper. It was needed when the necessary condition for prime $n$ that
$x^{3n}-ax^n-a\equiv 0 \pmod{n,x^3-ax-a}$ was satisfied. 
For example composite $n=13040299$ and $a=69121405197$.
We cannot determine
whether this g.c.d. is superfluous for the present cubic test.
Our testing includes it.

\section{Reinforcing testing over cubic polynomials}

Formally we compute $B\equiv x^{n-1} \pmod{n, f_a}$. The cases where $B\equiv 1$ are of no interest but
rather a non-trivial result of $B\equiv sx^2+tx+u$ is. We can then cheaply form and check the
stronger necessary condition: \[B^2+B+1 \equiv -x^2+x+a \pmod{n, x^3-ax-a}\,\,\, (**)\].
Given that $x$ and $B-1$ both have multiplicative inverses for prime $n$, the derivation of this working
\!\!\!\! $\pmod{n,f_a}$ throughout is as follows:
\[
\begin{array}{rcll}
x^{3}-ax- a&\equiv& 0&\textnormal{By definition from $(*)$.}\\
x^{3} &\equiv& ax+a&\textnormal{Add $ax+a$.}\\
x^{3} &\equiv& a(x+1)&\textnormal{Factor out $a$.}\\
x^{3n} &\equiv& a^n(x+1)^n&\textnormal{Raise to the $n^{th}$ power.}\\
x^{3n} &\equiv& a(x+1)^n&\textnormal{Fermat's Little Theorem on $a^n$.}\\
x^{3n} &\equiv& a(x^n+1)&\textnormal{Freshman's Dream for binomials.}\\
x^{3n}-x^3 &\equiv& a(x^n+1)-a(x+1)&\textnormal{Subtract identity.}\\
x^{3n}-x^3 &\equiv& a(x^n-x)&\textnormal{Collect terms.}\\
x^3(x^{3n-3}-1) &\equiv& ax(x^{n-1}-1)&\textnormal{Factor out $x$'s.}\\
x^3(B^3-1) &\equiv& ax(B-1)&\textnormal{Substitute to $B$.}\\
x^2(B^3-1) &\equiv& a(B-1)&\textnormal{Divide by $x$.}\\
x^2(B-1)(B^2+B+1) &\equiv& a(B-1)&\textnormal{Factor $B^3-1$.}\\
x^2(B^2+B+1) &\equiv& a&\textnormal{Divide by $B-1$.}\\
x^2(-x^2+x+a) &\equiv& a&\textnormal{Substitute the hypothetical identity from $(**)$.}\\
-x^4+x^3+ax^2 &\equiv& a&\textnormal{Multiply out.}\\
(-ax^2-ax)+(ax+a)+ax^2 &\equiv& a&\textnormal{Use the identity $x^3\equiv ax+a$.}\\
a &\equiv& a&\textnormal{Collect terms.}\\
\end{array}
\]

\section{Accelerating the search for a suitable B}

Recall the definition $B\equiv x^{n-1} \pmod{n, f_a}$.
Note that if $a=7+k(k-1)$ is prime then $k$ has to be $0$ or $1$ mod $3$.
More often than not, for prime $a$
and prime $n$ it is true that $B\equiv 1 \pmod{n, x^3-ax-a}$ implies $n^\frac{a-1}{3}\equiv 1 \pmod{a}$.
There is no known case of this implication for primes $a$ and $n$ to the contrary.
So the latter equation is a very useful screening test for the former equation.
This does not mean that composite $a$ cannot be utilised;
However in a practice using prime $a$ reduces the need to re-test with a different $a$ and potentially saves
a lot of computation of another $B$.
If $k\equiv 2 \pmod{3}$ is used then
$a$ is divisible by $3$.
If such an $a$ or another composite $a$ is used there is no need to check
$n^{a-1}\equiv 1 \pmod{a}$ or $n^\frac{a-1}{3}\equiv 1 \pmod{a}$.

\section{Cubic test summary}

The composite test for odd $n>1$ over $f_a=x^3-ax-a$ with prime $a=7+k(k-1)$ is as follows:
\begin{itemize}[itemsep=-1mm]
  {\item If $n$ is a perfect cube declare $n$ as a composite.}
  {\item If $a=n$ declare $n$ as a prime.}
  {\item If $n^{\frac{a-1}{3}}\equiv 1 \pmod{a}$ then start again with a new parametric prime $a$.}
  {\item Let $g=\gcd((2k-1)a(2a-1),n)$.}
  {\item If $g=n$ then start again with a new parametric prime $a$.}
  {\item If $1<g<n$ declare $n$ as a composite.}
  {\item Let $B\equiv x^{n-1} \pmod{n,f_a}$.}
  {\item If $B\equiv 1 \pmod{n,f_a}$ then start again with a new parametric prime $a$.}
  {\item If $B^2+B+1 \not\equiv -x^2+x+a \pmod{n,f_a}$ declare $n$ as a composite.}
  {\item Otherwise declare $n$ as a probable prime by the cubic test.}
\end{itemize}

\section{Verifications of the cubic test}

All of our verifications produced
{\em no pseudoprimes} for the cubic composite test.
The following verifications except the first two and the last
were run in the PARI/GP interpreter, some using Feitsma's base $2$ Fermat pseudoprimes $<2^{64}$
\cite{FEITSMA} and  Goutier's list of Carmichael numbers $<10^{22}$ \cite{GOUTIER}. 
See the Appendix I for the PARI/GP code. 
Appendix II contains references to The C++ Programming Language code.
\begin{enumerate}[label=(\roman*),itemsep=-1mm]
  {\item All prime and composite $a=7+k(k-1)$ for $a<n$ and odd $n<10^{8}$ (The C++ Programming Language).}
  {\item Minimal $k$ for odd $n<10^{15}$ (The C++ Programming Language).}
  {\item All $k\leq 600$ for Carmichael numbers $n<10^{22}$.}
  {\item First $20$ $B$: $B \not\equiv 1$ for base $2$ for the Fermat base 2 pseudoprimes $n<2^{64}$.}
  {\item All $k < n$ for Fermat base $2$ pseudoprimes $n<15\times 10^7$.}
  {\item All $a=2^r$ for all $r$ up to the multiplicative order of $2$ modulo $n$ for $n<10^{12}$.}
  {\item All $k<n=pq$ semi-primes: primes $p<76991$ and $q$ where $q=1+2j(p-1)$ for $4\leq j\leq 16$.}
  {\item All $k<n=pq$ semi-primes: primes $p<2729$ and $q$ where $q=1+2j(p^2-1)$ for $4\leq j\leq 16$.}
  {\item All $k<n=pq$ semi-primes: primes $p<2917$ and $q$ where $q=1+2j(p^2+p+1)$ for $4\leq j\leq 16$.}
  {\item First $20$ $B$: $B\not\equiv 1$ for odd $n<10^9$.}
  {\item 2.897,417,334,072 pseudorandom odd $n<2^{64}$ (The C++ Programming Language).}
\end{enumerate}

\section{Comparison with Baillie-PSW}
The $(n-1)^{th}$ power of $x$ is computed over $n$ and $f_a$. For small $a$ this
computation can be achieved using $O(6\log_2n)$ multiplications
and only $O(3\log_2n)$ modular reductions over $n$ in the main
with sundry multiplications by small numbers and additions.
With fast Fourier transform (FFT) arithmetic the test requires $O(3\log_2n)$ forward transforms
and $O(6\log_2n)$ inverse transforms.
This makes it very competitive against the Baillie-PSW tests \cite{BPSW}\cite{BPSW2},
which each require $O(4\log_2n)$ multiplications and
$O(4\log_2n)$ modular reductions; and $O(4\log_2n)$ forward FFT and $O(4\log_2n)$ inverse transforms.
\begin{figure}[h]
 \begin{center}
   \begin{tabular}{||l||c|c||c|c||}
      \hline
      Test & MUL & MOD & dFFT & iFFT\\
      \hline
      Cubic & 6 & 3 & 3 & 6\\
      BPSW & 4 & 4 & 4 & 4\\
      \hline
    \end{tabular}
    \caption{Test computational cost comparison}
 \end{center}
\end{figure}

Note that there would be easily found counterexamples to the Baillie-PSW tests
if the parameters were free and not just minimal, unlike the cubic composite test given in this paper.

\section{A quadratic test}

The only $a$ for which the cubic $f_a=x^3-ax-a$ is reducible is $a=8$.
 Then $f_8=(x+2)(x^2-2x-4)$.
This quadratic factor is irreducible and its discriminant is $20$.
For $\textnormal{jacobi} (5,n)=-1$ we can test base $-4$ Euler probable primality
in conjunction with $z^\frac{n+1}{2} \equiv \textnormal{jacobi} (-1,n) \pmod{n,z^2+3z+1}$.
As such it passes Feitsma's base $2$ Fermat pseudoprime list for $n<2^{64}$.
It is equivalent to Selfridge's \$500 challenge \cite{SELFRIDGE}
for a simultaneous base $2$ Fermat probable prime
and a pseudoprime with respect to the Fibonacci characteristic polynomial
$x^2-x-1$ for $n$ congruent to $2$ or $3$ modulo $5$.

\section{Conclusion}

The cubic test examined in this paper has not been shown to be a deterministic prime proving algorithm.
Moreover one may try to adapt the ideas presented in
Pomerance's paper \cite{POMERANCE} to show that we can expect to find pseudoprimes for
this cubic test. Like the Fermat probable prime test $b^{n-1}\equiv 1 \pmod{n}$
which results in non-trivial pseudoprimes
if the base $b$ is allowed to vary freely, tests based on Lucas sequences with respect to
$x^2-Px+Q$ are also weak in this sense; even if $P=c$, $Q=1$ and the Kronecker symbol of its discriminant $c^2-4$
over $n$ is $-1$.
However the cubic test over $x^3-ax-a$ presented here seems insusceptible to such a failing.

\section*{Appendix I -- PARI/GP code}
\footnotesize
Here is code for one way of computing $B=x^{n-1}\pmod{n,x^3-ax-a}$:
\begin{verbatim}
  {cubicB(n,k)=my(a=7+k*(k-1)%n,s=0,t=1,u=0,LEN=#binary(n-1));
    for(i=2,LEN,
      s2=a*sqr(s);t2=sqr(t);u2=sqr(u);
      st=2*a*s*t;tu=2*t*u;us=2*u*s;
      if(bittest(n-1,LEN-i),
        u=a*(s2+us+t2);s=s2+st+tu;t=u+u2+st,
        s=s2+us+t2;t=s2+st+tu;u=st+u2);
      s%=n;t%=n;u%=n);
    Mod(Mod(s*x^2+t*x+u,n),x^3-a*x-a);}
\end{verbatim}
\hrule
\vspace{6pt}

Here is code for the Cubic test using minimal $k$:
\begin{verbatim}
  {cubicTest(n)=my(k,a,B,g);
    if(n<2||ispower(n,3),
      return(0));
    k=0;B=1;
    while(B==1,
      k++;a=7+k*(k-1);
      while(!isprime(a)||Mod(n,a)^((a-1)/3)==1,
        k++;a=7+k*(k-1));
      if(a==n,
        return(1));
      g=gcd((2*k-1)*a*(2*a-1),n);
      if(1<g&&g<n,
        return(0));
      if(g==1,
        B=cubicB(n,k)));
    B^2+B+1==-x^2+x+a;}
\end{verbatim}
\hrule
\vspace{6pt}

{\bf Test (iii):}
\begin{verbatim}
  {V=readvec("~/Goutier/carm_10e22");}
  {for(k=1,600,
    a=7+k*(k-1);
    for(v=1,#V,
      n=V[v];
      if(gcd((2*k-1)*a*(2*a-1),n)==1,
        B=Mod(Mod(x,n),x^3-a*x-a)^(n-1);
        if(B^2+B+1==-x^2+x+a,
          print([n,k,a])))));}
\end{verbatim}
\hrule
\vspace{6pt}

{\bf Test (iv):}
\begin{verbatim}
  {V=readvec("~/Feitsma/PSP-2");}
  for(v=1,#V,n=V[v];
    k=0;a=7+k*(k-1);cnt=20;
    while(cnt,
      B=1;
      while(B==1,
        k++;a=7+k*(k-1);
        while(k%3==2||!ispseudoprime(a)||Mod(n,a)^((a-1)/3)==1||gcd((2*k-1)*a*(2*a-1),n)!=1,
          k++;a=7+k*(k-1));
        B=Mod(Mod(x,n),x^3-a*x-a)^(n-1));
      cnt--;
      if(B^2+B+1==-x^2+x+a,
        print([n,k,a])));}
\end{verbatim}
\hrule
\vspace{6pt}
\pagebreak
{\bf Test (v):}
\begin{verbatim}
  {V=readvec("~/Feitsma/PSP-2");}
  {for(v=1,#V,
    n=V[v];if(n>1.5*10^8,break};
    for(k=1,n/2,
      a=7+k*(k-1)%n;
      B=Mod(Mod(x,n),x^3-a*x-a)^(n-1);
      if(B^2+B+1==-x^2+x+a,
        print([n,k,a]))));}
\end{verbatim}
\hrule
\vspace{6pt}

{\bf Test (vi):}
\begin{verbatim}
  {V=readvec("~/Feitsma/PSP-2");}
  {for(v=1,#V,
    n=V[v];if(n>10^12,break};
    z=znorder(Mod(2,n));
    for(r=1,z,
      a=lift(Mod(2,n)^r);
      if(kronecker(4*a-27,n)==1,
        B=Mod(Mod(x,n),x^3-a*x-a)^(n-1);
        if(B^2+B+1==-x^2+x+a,
          print([n,a])))));}
\end{verbatim}
\hrule
\vspace{6pt}

{\bf Test (vii):}
\begin{verbatim}
  {tst(n,k)=my(a=7+k*(k-1),B=Mod(Mod(x,n),x^3-a*x-a)^(n-1));
  gcd((2*k-1)*a*(2*a-1),n)==1&&B^2+B+1==-x^2+x+a;}

  {tst1(p,q)=local(n=p*q,u=[],k,a,B);
  for(k=1,p,a=7+k*(k-1);B=Mod(Mod(x,p),x^3-a*x-a)^(n-1);
  if((n%(p-1)==1)||B^2+B+1==-x^2+x+a,u=concat(u,k)));Mod(u,p);}

  {tst2(p,q)=local(n=p*q,up,uq,k,V=[]);
  up=tst1(p,q);if(#up,uq=tst1(q,p);if(#uq,
  for(i=1,#up,for(j=1,#uq,k=lift(chinese(up[i],uq[j]));
  if(tst(n,k),V=concat(V,k))))));V=vecsort(V);
  if(#V,for(v=1,#V,t=V[v];print([n,t,7+t*(t-1)])));V;}

  {forprime(p=3,100000,for(k=4,16,q=1+2*k*(p-1);
  if(ispseudoprime(q),tst2(p,q))));
  print("\\\\ "round(gettime/1000)" seconds");}
\end{verbatim}
\hrule
\vspace{6pt}

{\bf Test (viii):}
\begin{verbatim}
  {forprime(p=3,100000,for(k=4,16,q=1+2*k*(p^2-1);
  if(ispseudoprime(q),tst2(p,q))));
  print("\\\\ "round(gettime/1000)" seconds");}
\end{verbatim}
\hrule
\vspace{6pt}

{\bf Test (ix):}
\begin{verbatim}
  {forprime(p=3,100000,for(k=4,16,q=1+2*k*(p^2+p+1);
  if(ispseudoprime(q),tst2(p,q))));
  print("\\\\ "round(gettime/1000)" seconds");}
\end{verbatim}
\hrule
\vspace{6pt}

{\bf Test (x):}
\begin{verbatim}
  {forcomposite(n=9,10^9,
    if(n%2==1&&!ispower(n,3),
      k=0;a=7+k*(k-1);cnt=20;
      while(cnt,
        B=1;
        while(B==1,
          k++;a=7+k*(k-1);
          while(k%3==2||!ispseudoprime(a)||Mod(n,a)^((a-1)/3)==1||gcd((2*k-1)*a*(2*a-1),n)!=1,
            k++;a=7+k*(k-1));
          B=Mod(Mod(x,n),x^3-a*x-a)^(n-1));
        cnt--;
        if(B^2+B+1==-x^2+x+a,
          print([n,k,a])))));}
\end{verbatim}
\hrule
\vspace{6pt}
\section*{Appendix II -- The C++ Programming Language code references}
\footnotesize
Gnu Multiprecision code for the cubic test:  \url{https://github.com/Boutoukoat/CubicPrimalityTest}
\vspace{6pt}
\hrule
\vspace{6pt}

{\bf Test (i)} and {\bf Test (ii)}  \url{https://github.com/Boutoukoat/CubicPrimalityTest}
\vspace{6pt}
\hrule
\vspace{6pt}

{\bf Test (xi)} Pseudorandom number generator used:
\begin{verbatim}
  uint64_t my_rand()
  {
      static uint64_t seed = 0x1234567812345678ull;
      seed = seed * 13 + 137; // linear congruential generator, period = 2^64
      return seed ^ (seed >> 13) ^ (seed << 13);   // bit shuffle
  }
\end{verbatim}
\end{document}